\documentclass[a4paper,10pt]{article}
\usepackage[utf8]{inputenc}
\usepackage{url}
\usepackage{fullpage}
\usepackage{authblk}

\newtheorem{theorem}{Theorem}

\newtheorem{proposition}[theorem]{Proposition}
\newtheorem{question}[theorem]{Question}

\usepackage{mathtools}

\DeclarePairedDelimiter\floor{\lfloor}{\rfloor}

\usepackage{mathrsfs}
\usepackage{amsfonts}
\usepackage{amssymb}
\usepackage{amsmath}
\usepackage{enumerate}

\title{A note on the simultaneous edge coloring \thanks{This work was supported by the ANR Project GrR (ANR-18-CE40-0032).}}
\author[1]{Nicolas Bousquet}
\author[2]{Bastien Durain}
\affil[1]{Univ. Grenoble Alpes, CNRS, Laboratoire G-SCOP, Grenoble-INP, Grenoble, France.\thanks{nicolas.bousquet@grenoble-inp.fr}}
\affil[2]{ENS Lyon, Département d'informatique, France. \thanks{bastien.durain@ens-lyon.fr}}

\begin{document}

\maketitle

\begin{abstract}
Let $G=(V,E)$ be a graph. A \emph{(proper) $k$-edge-coloring} is a coloring of the edges of $G$ such that any pair of edges sharing an endpoint receive distinct colors. A classical result of Vizing~\cite{Vizing64} ensures that any simple graph $G$ admits a $(\Delta(G)+1)$-edge coloring where $\Delta(G)$ denotes the maximum degreee of $G$. Recently, Cabello raised the following question: given two graphs $G_1,G_2$ of maximum degree $\Delta$ on the same set of vertices $V$, is it possible to edge-color their (edge) union with $\Delta+2$ colors in such a way the restriction of $G$ to respectively the edges of $G_1$ and the edges of $G_2$ are edge-colorings?
More generally, given $\ell$ graphs, how many colors do we need to color their union in such a way the restriction of the coloring to each graph is proper? 

In this short note, we prove that we can always color the union of the graphs $G_1,\ldots,G_\ell$ of maximum degree $\Delta$ with $\Omega(\sqrt{\ell} \cdot \Delta)$ colors and that there exist graphs for which this bound is tight up to a constant multiplicative factor. Moreover, for two graphs, we prove that at most $\frac 32 \Delta +4$ colors are enough which is, as far as we know, the best known upper bound.
\end{abstract}

\section{Introduction}
All along the paper, we only consider simple loopless graphs. In his seminal paper, Vizing proved in~g~\cite{Vizing64} that any simple graph $G$ can be properly edge-colored using $\Delta(G)+1$ colors (where $\Delta(G)$ denotes the maximum degreee of $G$). The \emph{union} of two graphs $G_1$ and $G_2$ on vertex set $V$ is the (simple) graph $G$ with vertex set $V$ and where $uv$ is an edge if and only if $uv$ is an edge of $G_1$ or an edge of $G_2$.
An edge coloring of $G$ is \emph{simultaneous with respect to $G_1$ and $G_2$} if its restrictions to the edge set of $G_1$ and to the edge set of $G_2$ are proper edge-colorings.
Recently, Cabello raised the following question\footnote{BIRS Workshop:  Geometric and Structural Graph Theory (2017). See e.g. \url{https://sites.google.com/site/sophiespirkl/birs-open-problems.pdf}}: given two graphs $G_1,G_2$ of maximum degree $\Delta$ on the same set of vertices $V$, does it always exist a simultaneous $(\Delta+2)$-edge coloring with respect to $G_1$ and $G_2$? Cabello proved that this property is satisfied if the intersection of $G_1$ and $G_2$ is regular~\cite{Cabello18}.
Using Vizing's theorem, one can easily notice that there exists a simultaneous $(2\Delta+1)$-edge coloring. From a lower bound perspective, no graph where $\Delta+2$ colors are needed is known.

Cabello introduced a generalization of this notion. Let $\ell$ graphs $G_1,G_2,\ldots,G_\ell$ and $G$ be their (edge) union. In other words, $uv$ is an edge of $G$ if and only if $uv$ is an edge of at least one graph $G_i$ with $i \le \ell$. An edge-coloring of $G$ is \emph{simultaneous with respect to $G_1,\ldots,G_\ell$} if its restriction to each graph $G_i$ is a proper edge-coloring. Cabello asked how many colors are needed to ensure the existence of a simultaneous coloring of $G$ with respect to each $G_i$. Let us denote by $\chi'(G_1,\ldots,G_\ell)$ the minimum number of colors needed to obtain a simultaneous coloring. And let $\chi'(\ell,\Delta)$ be the largest integer $k$ such that $k=\chi'(G_1,\ldots,G_\ell)$ for some graphs $G_1,\ldots,G_\ell$ of maximum degree (at most) $\Delta$. Vizing's theorem ensures that $\chi'(\ell,\Delta) \le \ell \Delta+1$ and Cabello exhibit a graph for which $\chi'(3,\Delta) \ge \Delta+5$ (with $\Delta=10$)~\cite{Cabello18}.
In this note, we prove that the order of magnitude of $\chi'(\ell,\Delta)$ is $\Theta(\sqrt{\ell} \Delta)$. More precisely, we prove that the following statement holds:

\begin{theorem}\label{thm:sqrt}
 \[ \chi'(\ell,\Delta) \le 2 \sqrt{2 \ell} \Delta - \sqrt{2 \ell} + 2. \]
\end{theorem}

We claim that this bound is tight up to a constant multiplicative factor. Let $\ell \in \mathbb{N}$ and $\Delta$ be an even value. Let $G:=S_{1,k \Delta}$ be the star with $k\Delta$ leaves, where $k= \lfloor \sqrt{ \frac{\ell}{2} } \rfloor$. Partition the edges of $G$ into $2k$ sets $A_1,\ldots,A_{2k}$ of size $\frac \Delta 2$. 
For every pair $i,j$, create the graph $G_{i,j}$ with edge set $A_i \cup A_j$. Note that each graph $G_{i,j}$ has maximum degree $\Delta$ since by construction the set of edges $A_i$ induces a graph of maximum degree $\Delta/2$. Moreover the total number of graphs $G_{i,j}$ is $2k(2k-1)/2 \le \ell$. Finally, by construction, every pair of edges of $G$ appears in at least one graph $G_{i,j}$. So in order to obtain a simultaneous coloring, we need to color all the edges of $G$ with different colors (since all the edges are incident to the center of the star). So:
\begin{proposition}\label{prop:lb}
 \[ \chi'(\ell,\Delta) \ge \floor*{ \sqrt{\frac \ell 2} } \Delta.\] 
\end{proposition}
Note that, for $\ell=3$ and the graph $S_{1,3\floor*{ \Delta/2 }}$, a similar construction ensures that $\chi'(3,\Delta) \ge 3 \floor*{ \frac \Delta 2 }$, improving the lower bound of $\Delta+5$. Indeed, let us partition the edges of the star into three sets $A_1,A_2,A_3$ of size $\floor*{ \Delta/2 }$. We similarly define for every $i\neq j$ the graph $G_{i,j}$ with edge set $A_i \cup A_j$. Each graph $G_{i,j}$ has maximum degree $\Delta$ and every pair of edges appears in at least one graph $G_{i,j}$. So an edge coloring of $S_{1,3\floor*{ \Delta/2 }}$ simultaneous with respect to $G_{1,2},G_{1,3}$ and $G_{2,3}$ is a proper edge coloring of $S_{1,3\floor*{ \Delta/2 }}$.

When $\ell=2$, a careful reading of the proof of Theorem~\ref{thm:sqrt} permits to remark that we can improve the trivial $(2\Delta+1)$ upper bound into $2\Delta$.
We prove the following much better upper bound with a different technique:

\begin{theorem}\label{thm:linear}
  \[ \chi'(2,\Delta) \le \floor*{ \frac{3}{2}  \Delta +4 }. \]
\end{theorem}
As far as we know, it is the best known upper bound.

\section{Proof of Theorem~\ref{thm:sqrt}}

Let $G_1,\ldots,G_\ell$ be $\ell$ graphs of maximum degree $\Delta$. Let us partition the set of edges of $G = \cup_{i=1}^\ell G_i$ into two sets (all along the paper, the notation $\cup$ stands for edge union). The \emph{multiplicity} of an edge $e$ is the number of graphs $G_i$ with $i \le \ell$ on which $e$ appears.
For some fixed $k$, the set $E_1$ is the set of edges with multiplicity at least $k$ and $E_2$ is the set of edges with multiplicity less than $k$. We will optimize the value of $k$ later. (Note that we do not necessarily assume that $k$ is an integer). For every $i \in \{1,2\}$, let us denote by $H_i$ the graph $G$ restricted to the edges of $E_i$. Note that $G= H_1 \cup H_2$. 

We claim that the graph $H_1$ has degree at most $\ell \Delta / k $. Indeed, let $u$ be a vertex and $E_1(u)$ be the set of edges of $H_1$ incident to it. Since every edge of $H_1$ has multiplicity at least $k$ and at most $\ell \Delta$ edges (with multiplicity) are incident to $u$ in $G$, at most $\frac \ell k \Delta$ different edges are in $E_1(u)$. So $H_1$ has maximum degree $\frac \ell k \Delta$. By Vizing's theorem, $H_1$ can be properly edge-colored with $ (\frac \ell k \Delta+1)$ colors. 

Let us now prove that $H_2$ can be simultaneously edge-colored with $2 k (\Delta-1)+1$ colors. Let us prove it by induction on the number of edges of $H_2$. The empty graph can indeed be edge-colored with $2 k (\Delta-1)+1$ colors. Let $e=uv$ be an edge of $H_2$. By induction, there exists a simultaneous coloring $c'$ of $H_2 \setminus e$ with $2k (\Delta-1)+1$ colors. Let us prove that $c'$ can be extended into a simultaneous coloring of $H_2$. Without loss of generality, we can assume that $e$ is an edge of the graphs $G_1,\ldots,G_{r}$ with $r < k$ and is not an edge of $G_{r+1},\ldots,G_\ell$. Let $F$ be the set of edges of $G_1,\ldots,G_r$ incident to $u$ or to $v$ distinct from $e$. By assumption, there are at most $2 r (\Delta-1)$ such edges ($2(\Delta-1)$ in each graph). Since $r < k$, at most $2 k (\Delta-1)$ edges are in $F$. So there exists a color $a$ that does not appear in $F$. The edge $e$ can be colored with $a$ without violating any constraints. 
It holds by choice of $a$ for $G_i$ with $i\le r$ and it holds since $e \notin G_i$ for $i > r$.
 
 So $\chi'(\ell,\Delta) \le \frac \ell k \Delta + 2  k  (\Delta-1)+2$ colors. We finally optimize the integer $k$ which minimize the number of colors. We want to minimize $\frac \ell k + 2k$ which is minimal when $k= \sqrt{\frac \ell 2}$. It finally ensures that $\chi'(\ell,\Delta) \le 2 \sqrt{2 \ell} \Delta - \sqrt{2 \ell} + 2$, which completes the proof of Theorem~\ref{thm:sqrt}.

 \section{Proof of Theorem~\ref{thm:linear}}

 Let $G_1,G_2$ be two graphs of maximum degree $\Delta$ and let $G$ be their union. Let $E_2$ be the  edges that appear in both graphs and, for every $i \in \{ 1,2\}$, let $E_1^i$ be the set of edges that appears only in $G_i$. Let us denote by $H_2$ (resp. $H_1^i$) the graph restricted to the edges of $E_2$ (resp. $E_1^i$).

 
 
 
 For every vertex $v$ and every graph $H$, we denote by $degH(v)$ the degree of $v$ in $H$.
 Let $H$ be a graph and $f,g$ be two functions from $V(H)$ to $\mathbb{R}^+$.  A \emph{$(g,f)$-factor} of $H$ is an edge-subgraph $H'$ of $H$ such that every vertex $v$ satisfies $g(v) \le deg_{H'}(v) \le f(v)$ .
 Kano and Saito proved in~\cite{KanoS83} that the graph $H$ admits a $(g,f)$-factor if 
 \begin{enumerate}[(i)]
  \item $f$ and $g$ are two integer valued functions, and
  \item for every vertex $v$, $g(v) < f(v)$, and
  \item there exists  a real number $\theta$ such that $0 \le \theta \le 1$ and for every vertex $v$, $g(v) \le \theta \cdot deg_{H}(v) \le f(v)$.
 \end{enumerate}

 Let $1 \le i \le 2$. We will extract from $H_1^i$ a $(g,f)$-factor where $g(v)=\lceil \frac{deg_{H_1^i}(v)}{2} -1\rceil$ and $f(v)=\lceil \frac{deg_{H_1^i}(v)}{2} \rceil$. The points (i) and (ii) are satisfied. Moreover, by choosing $\theta = \frac 12$, (iii) is also satisfied. Thus by~\cite{KanoS83}, the graphs $H_1^1$ and $H_1^2$ admit $(g,f)$-factors.
 For $i \le 2$, let $K_1^i$ be a $(g,f)$-factor of $H_1^i$. 
 For every $i$, let $L_1^i = H_1^i \setminus K_1^i$. Let $L= L_1^1 \cup L_1^2$ and $R = H_2 \cup K_1^1 \cup K_1^2$. Note that $G= L \cup R$. Let us now color these two graphs. 
 \smallskip
 
Let us first prove that $L$ can be colored with $\floor*{ \frac{\Delta}{2} }+2$ colors. For every $i$, the graph $L_1^i$ has maximum degree at most $\floor*{ (\Delta/2)+1 }$ since every vertex $v$ of $K_1^i$ has degree at least $\lceil (deg_{H_1^i}(v)/2)-1  \rceil$. By Vizing's theorem, the graph $L_1^i$ can be colored with at most $\floor*{ \frac{\Delta}{2}+1 } +1 = \floor*{ \frac{\Delta}{2} } +2$ colors. Since the edges of $L_1^1$ and $L_1^2$ are disjoint, $L=L_1^1 \cup L_1^2$ can be simulteaneously colored with $\floor*{ \frac{\Delta}{2} } +2$ colors (the same set of colors can be re-used for each graph).
\smallskip

Let us now color the graph $R$. 
Let $v$ be a vertex of $R$. Let us denote by $d$ the degree of $v$ in $H_2$. Since edges of $H_2$ are in both $G_1$ and $G_2$, the vertex $v$ has degree at most $\Delta-d$ in both graphs $H_1^1$ and $H_1^2$. Since the graphs $K_1^1$ and $K_1^2$ are $(g,f)$-factors of respectively $H_1^1$ and $H_1^2$, the degree of the vertex $v$ is at most $\lceil \frac{\Delta-d}{2} \rceil$ in each graph. So the degree of $v$ in the graph $R$ is at most $d +2 \lceil \frac{\Delta-d}{2} \rceil \le \Delta+1$. By Vizing's theorem, the graph $R$ can be colored using at most $\Delta+2$ colors.
\smallskip

Since $G=L \cup R$, we can find a simultaneous edge-coloring with respect to $G_1$ and $G_2$ using at most $\floor*{ \frac{\Delta}{2} } +2 + \Delta+2 = \floor*{ \frac 32 \Delta }+4$ colors.

\section{Conclusion}
Theorem~\ref{thm:sqrt} and Proposition~\ref{prop:lb} ensures that the following holds:
\[ \floor*{ \sqrt{\frac \ell 2} } \Delta \le \chi'(\ell,\Delta) \le 2 \sqrt{2 \ell} \Delta - \sqrt{2 \ell} + 2. \]
Closing the multiplicative gap of $4$ between lower and upper bound is an interesting open problem. For $\ell=2$, we still do not know any graph for which $\chi'(2,\Delta) > \Delta+1$. Cabello asked the following question that is still widely open despite the progress obtained in Theorem~\ref{thm:linear}:
\begin{question}[Cabello]
Is it true that
\[ \chi'(2,\Delta) \le \Delta+2 \ \ ? \]
\end{question}

 \paragraph{Acknowledgements.} The authors want to thank Louis Esperet for fruitful discussions and suggestions.

 \bibliographystyle{abbrv}

\end{document}